\documentclass[12pt]{amsart}
\usepackage{amsmath,amssymb,latexsym}
\usepackage[a4paper,total={16cm,23.5cm},top=3cm,left=2.5cm]{geometry}

\theoremstyle{plain}
\newtheorem{theorem}{Theorem}[section]
\newtheorem{proposition}[theorem]{Proposition}
\newtheorem{fact}[theorem]{Fact}
\newtheorem{lemma}[theorem]{Lemma}
\newtheorem{corollary}[theorem]{Corollary}

\theoremstyle{definition}
\newtheorem{definition}[theorem]{Definition}
\newtheorem{remark}[theorem]{Remark}
\newtheorem{expl}[theorem]{Example}
\newtheorem{frage}[theorem]{Question}

\def\bsp{\begin{expl}}
\def\ebsp{\end{expl}}
\def\beh{\begin{claim}}
\def\ebeh{\end{claim}}
\def\defn{\begin{definition}}
\def\edefn{\end{definition}}
\def\satz{\begin{theorem}}
\def\esatz{\end{theorem}}
\def\tats{\begin{fact}}
\def\etats{\end{fact}}
\def\kor{\begin{corollary}}
\def\ekor{\end{corollary}}
\def\lmm{\begin{lemma}}
\def\elmm{\end{lemma}}
\def\bem{\begin{remark}}
\def\ebem{\end{remark}}
\def\bew{\begin{proof}}
\def\ebew{\end{proof}}
\def\satzli{\begin{proposition}}
\def\esatzli{\end{proposition}}
\def\frg{\begin{frage}}
\def\efrg{\end{frage}}

\def\Ind#1#2{#1\setbox0=\hbox{$#1x$}\kern\wd0\hbox to 0pt{\hss$#1\mid$\hss}
\lower.9\ht0\hbox to 0pt{\hss$#1\smile$\hss}\kern\wd0}
\def\ind{\mathop{\mathpalette\Ind{}}}
\def\Notind#1#2{#1\setbox0=\hbox{$#1x$}\kern\wd0\hbox to 0pt{\mathchardef
\nn="3236\hss$#1\nn$\kern1.4\wd0\hss}\hbox to 0pt{\hss$#1\mid$\hss}\lower.9\ht0
\hbox to 0pt{\hss$#1\smile$\hss}\kern\wd0}

\def\F{\mathcal F}
\def\B{\mathcal B}

\begin{document}
\title{Bad Groups}
\author{Frank O. Wagner}
\address{Universit\'e de Lyon; Universit\'e Claude Bernard Lyon 1; CNRS; Institut Camille Jordan UMR5208, 43 bd du 11 novembre 1918, 69622 Villeurbanne Cedex, France}

\email{wagner@math.univ-lyon1.fr}

\date{\today}
\thanks{Partially supported by ANR-13-BS01-0006 ValCoMo}
\begin{abstract}There is no bad group of Morley rank $2n+1$ with an abelian Borel subgroup of Morley rank $n$. In particular, there is no bad group of Morley rank $3$ (O. Fr\'econ).\end{abstract}
\maketitle

\section{Introduction}
The algebraicity conjecture, proposed independently by Gregory Cherlin \cite{Ch79} and Boris Zilber \cite{Zi77} forty years ago, states that a simple group of finite Morley rank should be an algebraic group over an algebraically closed field. Despite much effort, it is still open. In fact, it is currently believed that in the so-called degenerate case, where there is no involution, a counter-example may exist. It has been shown ten years ago by Alt\i nel, Borovik and Cherlin \cite{ABC08} that if there is an infinite elementary abelian $2$-subgroup, the conjecture holds. However, the inductive approach for the remaining case where a connected $2$-Sylow subgroup is abelian and divisible got stuck at the initialisation stage: it is very difficult to eliminate small non-algebraic configurations.

In the paper \cite{Ch79}, Cherlin classified connected groups of small Morley rank: In rank $1$ they are abelian after a result of Reinecke \cite{Re75}, in rank $2$ they are soluble, and in rank $3$ they are either soluble, PSL$_2(K)$ for some interpretable algebraicaly closed field $K$, or what he called a {\em bad group}. Bad groups were further studied by Nesin \cite{Ne89}, Corredor \cite{Co89} and Borovik-Poizat \cite{BP90}. A connected group of finite Morley rank is {\em bad} if all its soluble connected subgroups are nilpotent. In a minimal bad group, the maximal nilpotent connected subgroups (called {\em Borel subgroups}) are definable, self-normalizing, conjugate, and the conjugates cover the whole group; moreover the group does not have an involution. It is easy to see \cite[p.\ 91]{Po87} that the Morley rank of a Borel subgroup is strictly less than half the rank of the ambient bad group, but essentially little progress has been made for more than 25 years until 2016, when Olivier Fr\'econ \cite{Fr16} proved the non-existence of bad groups of Morley rank three, and thus completed Cherlin's classification from almost 40 years ago.

Fr\'econ's proof is restricted to rank $3$, and our generalisation to rank $2\,RM(B)+1$ (where $B$ is any abelian Borel subgroup) is absolutely minor. However, given that bad groups are one of the main obstacles to a full proof of the algebraicity conjecture (at least in the presence of involutions), it seems important to analyze the scope and current limits of Fr\'econ's arguments. This note hopes to contribute to this study.

\section{Full Frobenius groups}
Recall that a subgroup $B$ of a group $G$ is {\em malnormal} if $B\cap B^g=\{1\}$ for any $g\in G\setminus B$. This means that $B$ is self-normalizing and intersects all its distinct conjugates trivially. Thus, for any $b\in B\setminus\{1\}$, if $b^g\in B$ for some $g\in G$, then $g\in B$; in particular $C_G(b)\le B$.
A malnormal subgroup is often called a {\em Frobenius complement} and the ambient group a {\em Frobenius group}; Frobenius' Theorem states that in a finite group, a malnormal subgroup $B$ has a normal complement, called the {\em Frobenius kernel} and consisting of all elements outside all conjugates of $B\setminus\{1\}$. It is easy to deduce that Frobenius' Theorem also holds for locally finite groups, and by the transfer principle for all algebraic groups. The following definition, due to Jaligot \cite{Ja01}, isolates a particularly pathological configuration.
\defn A Frobenius group is {\em full} if the conjugates of the Frobenius complement cover the whole group.
A malnormal subgroup whose conjugates cover the whole group is a {\em full Frobenius complement}.
\edefn
We shall see that a connected full Frobenius group of finite Morley rank has a simple definable subgroup which is still a full Frobenius group, and hence contradicts the algebraicity conjecture. One would thus like to show that they do not exist.

For the rest of this section $G$ will be a full Frobenius group of finite Morley rank, and $B$ a malnormal subgroup. We shall call the conjugates of $B$ the Borel subgroups of $G$.

\lmm\label{involution} $G$ has no involutions. Thus, every element $g$ of $G$ has a unique square root, which is contained in every definable subgroup containing $g$. If $g\not=1$, it is not conjugate to its inverse $g^{-1}$. If $G$ is connected, then $C_G(g)$ is infinite for all $g\in G$\elmm
\bew If $i$ and $j$ are two involutions in different conjugates of $B$, then both invert the element $ij$, and hence normalize the unique Borel $B'$ containing $ij$. But then both $i$ and $j$ are contained in $B'$, a contradiction.

Let $H$ be the smallest definable subgroup containing $g$. Then $H\le Z(C_G(g))$ is abelian, and $x\mapsto x^2$ is an injective homomorphism. So it must be surjective, since its image has the same Morley rank and degree as $H$. Hence $g$ has a unique square root $h$ in $H$. Every other square root $h'$ of $g$ must commute with $g$ and normalize $H$, and hence commute with the unique square root $h$ of $g$ in $H$. But then $(h'h^{-1})^2=1$, and $h=h'$.

If $g^h=g^{-1}$, then $h^2\in C_G(g)$, so $h\in C_G(g)$, and $g=g^{-1}$. Since there are no involutions, $g=1$.

Finally, if $g$ has a finite centraliser, then the conjugacy class $g^G$ is generic, as is $(g^{-1})^G$; if $G$ is connected then $g$ is conjugate to its inverse, a contradiction.\ebew

Note that all the assertions of Lemma \ref{involution} follow from the first one, and thus hold in any group of finite Morley rank without involutions.

\lmm\label{finite} Every finite or soluble subgroup is contained in a Borel. In particular, $G$ has no finite or soluble normal subgroup.\elmm
\bew It is clear that a Borel $B$ cannot contain a normal subgroup, as the normalizer of any subgroup of $B$ must be contained in $B$ by malnormality.

Let $F$ be a finite subgroup not contained in a Borel; we may assume that is it minimal such. If $B$ is a Borel intersecting $F$ non-trivially, then $B\cap F$ is malnormal in $F$, so $F$ is a finite Frobenius group. By minimality, its kernel $N$ is contained in a Borel subgroup $B'$. As $N$ is non-trivial, $F\le N_G(N)\le B'$, a contradiction.

Let $S$ be a non-trivial soluble subgroup. Then $S$ has a non-trivial abelian normal subgroup $A$, which must be contained in a Borel $B$, and $S\le N_G(A)\le B$.\ebew

\lmm $G$ has a unique non-trivial minimal normal definable subgroup, which is simple.\elmm
\bew As $G$ has no soluble normal subgroup by Lemma \ref{finite}, its socle is a finite product of definable simple groups $S_1\times\cdots\times S_n$ by \cite[p.\ 97]{Po87}. Consider $s\in S_1\setminus\{1\}$ and a Borel $B$ containing $s$. Then $S_i\le C_G(s)\le B$ for any $i>1$. But if $s'\in S_2\setminus\{1\}$, then $S_1\le C_G(s')\le B$ as well, and $B$ contains a non-trivial normal subgroup, contradicting malnormality. It follows that $n=1$.\ebew

\lmm\label{full}\begin{enumerate}
\item $B$ is infinite.
\item If $G$ is connected, so is $B$.
\item If $H$ is a connected definable subgroup of $G$ which is not contained in a Borel, then $H$ is a full Frobenius group, and $B\cap H$ is a full Frobenius complement in $H$, for any Borel $B$ of $G$ with $B\cap H\not=\{1\}$.
\item If $H$ is a connected definable subgroup containing a Borel $B$, then $H$ is a full Frobenius complement for $G$.
\item If $G$ is connected, $B'$ is another full Frobenius complement and $B\cap B'$ is nontrivial, then $B\cap B'$ is again a full Frobenius complement. In particular there is a unique minimal full Frobenius complement.
\item If $N$ is a normal subgroup of $G$, then $G=NB$.\end{enumerate}\elmm
\bew\begin{enumerate}
\item Let $g\in G^0$. Then $C_{G^0}(g)$ is infinite, and contained in the Borel of $g$.
\item If $B$ were not connected, then the union of the conjugates of $B^0$ and the union of the conjugates of $B\setminus B^0$ would be two disjoint generic subsets of $G$.
\item Let $(B_i:i\in I)$ be the non-trivial intersections with $H$ of the conjugates of $B$. Then $H$ is the disjoint union of the $B_i$, which are malnormal in $H$.
The union of the $H$-conjugates of $B_i$ is generic in $H$ for all $i\in I$; by connectedness the $B_i$ are all $H$-conjugate.
\item If $g\in G$ and $h\in H\setminus\{1\}$ with $h^g\in H$, then there are two Borel subgroups $B'$ and $B''$ of $H$ with $h\in B'$ and $h^g\in B''$. But there is $h'\in H$ with $B'^{h'}=B''$, so $h^g\in B'^g\cap B'^{h'}$. It follows that $B'^g=B'^{h'}$, and $gh'^{-1}\in N_G(B')=B'\le H$. Hence $H$ is malnormal in $G$; it is clear that its conjugates cover $G$.
\item We may assume $B\not\le B'$. Then $B$ is connected by (2), and $B\cap B'$ is a full Frobenius complement in $B$ by (3). The result now follows from fullness of $B$.
\item $N$ is infinite by Lemma \ref{finite}. But then $N^0$ cannot be contained in a Borel by normality, and must be a full Frobenius group itself; if $B$ is a Borel of $G$ and $B\cap N^0$ is a Borel of $N^0$, then for any $g\in G$ there is $h\in N^0$ with $(B\cap N^0)^g=(B\cap N^0)^h$, whence 
$gh^{-1}\in N_G(B)=B$. Thus $G=N^0B=NB$.\qedhere\end{enumerate}\ebew

If in (3) we choose $H$ of minimal Morley rank not contained in some Borel subgroup, we obtain a full Frobenius group whose connected proper definable subgroups are all contained in some Borel of $H$. In fact all definable subgroups are contained in some Borel of $H$: This is clear for finite groups by Lemma \ref{finite}, and for a proper infinite definable subgroup $K$ there is a Borel $B$ of $H$ with $K^0\le B$, and $K\le N_G(K^0)\le B$ by malnormality. In particular the Borel subgroups of $H$ are precisely its maximal definable proper subgroups, which are connected by (2).

\defn A full Frobenius group is {\em special} if its Borel subgroups are precisely its maximal proper definable subgroups.\edefn
In particular a special full Frobenius group $G$ must be simple, whence connected. If moreover its Borel subgroups are nilpotent, then $G$ is a bad group in the sense of the introduction. Conversely, a minimal bad group gives rise to a special full Frobenius group with nilpotent complement \cite{Ne89,Co89,BP90}. Thus if there are no special full Frobenius groups, neither are there full Frobenius groups nor bad groups.

\section{Involutive Automorphisms}
In this section we shall study a group $G$ of finite Morley rank and without involutions, together with an involutive automorphism $\sigma$ of $G$. Let $F$ be the set of elements invariant under $\sigma$, and $I$ the set of elements inverted by $\sigma$.
\lmm\label{l5}\begin{enumerate}
\item No non-trivial lement of $F$ is conjugate to an element of $I$.
\item For every $g\in G$ there are unique $f\in F$ and $i\in I$ with $g=fi$. In particuliar, if $G$ is connected, so is $F$.
\end{enumerate}\elmm
\bew\begin{enumerate}\item The sets $\{g\in G:\sigma(g)\mbox{ is conjugate to }g\}$ and $\{g\in G:\sigma(g)\mbox{ is conjugate to }g^{-1}\}$ are $G$-invariant and contain $F$ and $I$, respectively. They must have empty intersection, as no element is conjugate to its inverse.
\item Put $g(x) = x^{-1}\sigma(x)$. Then $g(x)\in I$ for all $g\in G$. So the unique square root $z$ of $g(x)^{-1}$ must also be in $I$. Then
$$g(x)=x^{-1}\sigma(x)=z^{-2}=z^{-1}\sigma(z) = g(z),$$
so $xz^{-1}=\sigma(xz^{-1})\in F$. Hence $G=FI$.

If $fi=f'i'$ with $f,f'\in F$ and $i,i'\in I$, then $i'i^{-1}\in F$, and
$$i'i^{-1}=\sigma(i'i^{-1})=i'^{-1}i,$$
whence  $i'^2=i^2$ and $i'=i$. It follows that $f=f'$.

Finally, by uniqueness of the decomposition, the generic types of $G$ are in bijection with the independent products of the types of maximal rank in $F$ and in $I$. If $G$ is connected, $F$ and $I$ must have a unique type of maximal rank; in particular $F$ is connected.\qedhere\end{enumerate}\ebew
Some special cases of the following Proposition have been shown in \cite[p.\ 393]{BN94} and \cite[p.\ 128]{Ja01}.
\satzli\label{p6} A special full Frobenius group has no non-trivial definable involutive automorphism. If $\sigma$ is a definable involutive automorphism of a connected full Frobenius group $G$, then $RM(F\cap B)<RM(B)$ for any minimal full complement $B$, and $RM(B)\ge2$. Moreover, $$RM(I)\le\frac{RM(G)+RM(B)}2-1<\frac34\,RM(G)-1.$$\esatzli
\bew Let $\sigma$ be a definable non-trivial involutive automorphism. Then $F$ is a proper definable subgroup of $G$.

If $G$ is special, $F$ is contained in a single Borel $B$. As $G=FI$ there must be some $i\notin B$ inverted by $\sigma$; if $B'$ is the Borel containing $i$, then $\sigma(B)$ is again a maximal proper subgroup, and hence conjugate to $B'$. As $\sigma(B')\cap B'$ is non-trivial, $\sigma(B')=B'$ and $\sigma$ stabilises $B'$. Thus $\sigma$ is a definable involutive automorphism of $B'$ without fixed points; by Lemma \ref{l5}(2) it inverts $B'$. But $B'$ is conjugate to $B$, contradicting Lemma \ref{l5}(1). The same argument works if $G$ is connected and $F$ is contained in a single minimal full complement, noting that $\sigma(B')$ is again a minimal full complement and hence conjugate to $B'$.

So $F$ is a full Frobenius group with full complement $F\cap B$ by Lemma \ref{full}(3), where $B$ is a minimal full complement intersecting $F$ non-trivially. No Borel which does not intersect $F$ can contain a point of $I$, as otherwise it would itself be inverted by $\sigma$, again contradicting Lemma \ref{l5}. Moreover $F\not\le B$ and both are connected, so $1\le RM(F\cap B)<RM(B)$~; in particular $RM(B)\ge2$.

If $\B$ is the family of conjugates of $B$ intersecting $F$ non-trivially, then 
$$RM(\B)=RM(F)-RM(F\cap B),$$
and for $B'\in\B$ we have $RM(I\cap B')=RM(B)-RM(F\cap B)$. Hence 
$$\begin{aligned}RM(I)&=RM(\B)+RM(I\cap B)\\
&=RM(F)-RM(F\cap B)+RM(B)-RM(F\cap B)\\
&=RM(F)+RM(B)-2\,RM(F\cap B)\,;\\
2\,RM(I)&=RM(I)+RM(F)+RM(B)-2\,RM(F\cap B)\\
&=RM(G)+RM(B)-2\,RM(F\cap B).\end{aligned}$$
Therefore 
$$\begin{aligned}RM(I)&=\frac{RM(G)+RM(B)}2-RM(F\cap B)\\
&\le \frac{RM(G)+RM(B)}2-1.\end{aligned}$$
Now $RM(G)>2\,RM(B)$, as otherwise for $g\notin B$ the double coset $BgB$ would be generic in $G$ and equal to $Bg^{-1}B$, which would imply the presence of involutions (see \cite[p.\ 91]{Po87}). Hence $$RM(I)<\frac34\,RM(G)-1.\qedhere$$
\ebew

\section{Almost twistedly normal subsets}
In this and the following section we shall consider particular subsets of a group $G$ which produce involutive automorphisms, and hence cannot exist in a special full Frobenius group.
\defn Two definable sets $X$ and $Y$ are {\em almost equal}, denoted $X\sim Y$, if $RM(X\triangle Y)< RM(X)$.\edefn
Note that $\sim$ is an equivalence relation, and $X\sim Y$ implies $RM(X)=RM(Y)=RM(X\cap Y)$. 
\defn A definable subset $X$ of an $\omega$-stable group $G$ is {\em twistedly normal} if for all $g\in G$ there is $h\in G$ with $gX=Xh$. The subset $X$ is {\em almost twistedly normal} if for all $g\in G$ there is $h\in G$ with $gX\sim Xh$.\edefn
We shall first show that an almost twistedly normal subset of a group gives rise to a twistedly normal subset.

\lmm\label{lsch} Let $G$ be a group acting definably on a set $X$ in an $\omega$-stable structure. Let $Y$ be a definable subset of $X$ such that $gY\sim Y$ for all $g\in G$. Then there is a $G$-invariant $Z\subseteq X$ with $Z\sim Y$.\elmm
\bew Let $p_1,\ldots,p_n$ be the generic types of $G$. We add the parameters necessary to define $G$, $X$, $Y$ and the $p_i$ to the language, and put
$$Z_i=\{x\in X:x\in gY\mbox{ for some/every }g\models p_i\mbox{ with }g\ind x\}.$$
Then $Z_i$ is $\emptyset$-definable by definability of types. Let $x\in Y$ be of maximal rank, and $g\models p_i$ generic over $x$. As $RM(Y\triangle gY)<RM(Y)$, we have $x\in gY$ and $x\in Z_i$. In particular $RM(Z_i)\ge RM(Y)$. Conversely, if $x\in Z_i$ is of maximal rank, there is $g\models p_i$ independent of $x$ with $x\in gY$. Hence
$$RM(Y)\le RM(Z_i)=RM(x)=RM(x/g)\le RM(gY)=RM(Y).$$
Thus $x$ is of maximal rank in $gY$, and $x\in gY\cap Y$, whence $x\in Y$. So $Y\sim Z_i$.

Finally, put $Z=\bigcap_iZ_i$. Clearly $Z$ is almost equal to $Y$, and $Z$ is invariant under all generics of $G$, and hence $G$-invariant.\ebew
\satzli Let $X$ be an almost twistedly normal subset of $G$. Then there is a twistedly normal subset $Y$ with $X\sim Y$.\esatzli
\bew We consider the action of $G\times G$ on $G$ given by $(g,h)(x)=gxh^{-1}$. By definability of rank, the subgroup
$$H=\{(g,h)\in G\times G:gXh^{-1}\sim X\}$$
is definable. By Lemma \ref{lsch} there is an $H$-invariant subset $Z$ of $G$ with $X\sim Y$; clearly $Y$ is twistedly normal.\ebew
\satzli\label{p7}\cite{Fr16} Let $G$ be a simple group of finite Morely rank, and $X$ an infinite non-generic definable subset of $G$. If $X$ is almost twistedly normal, then there is a unique definable automorphism $\sigma$ of $G$ such that $gX\sim X\sigma(g)$ for all $g\in G$.\esatzli
\bew Put $F=\{g\in G:gX\sim X\}$, a definable subgroup of $G$. As $\emptyset\not\sim X\not\sim G$ by hypothesis, $F$ is a proper subgroup of $G$. If $x\in F$ and $g\in G$, choose $h\in G$ with $gX\sim Xh$. Then
$$xgX\sim xXh\sim Xh\sim gX.$$ 
Hence $x^g\in F$ and $F$ is normal in $G$, whence trivial by simplicity.

By a similar argument, $F'=\{g\in G:X\sim Xg\}$ is normalised by all $h\in G$ such that there is some $g\in G$ with $gX\sim Xh$. So there is an injective definable homomorphism from $G$ to $N_G(F')/F'$ such that $gX\sim Xh$ for all $h\in\sigma(g)$. But 
$$RM(G)\ge RM(N_G(F'))\ge RM(N_G(F')/F')\ge RM(G)$$
by injectivity, and equality holds all the way. As $G$ is connected, we have $G=N_G(F')$, and $F'$ must be trivial as well by simplicity of $G$.\ebew

\section{Fr\'econ Elements}
In this section we shall consider a connected full Frobenius group $G$ with full complement $B$.
Following Fr\'econ, we shall call call a double translate $gBh^{-1}$ a {\em line}. Note that any line is a left (or right) translate of a conjugate of $B$, and that any two distinct lines intersect in at most one point. Given any two distinct points $x$ and $y$, there is a unique line $xB'$ containing them both, where $B'$ is the conjugate of $B$ containing $x^{-1}y$. 

If $X$ is a definable subset, we shall say that a line $d$ is {\em generically contained} in $X$ if $RM(d\setminus X)<RM(d)$. As $B$ and $d$ have Morley degree $1$, this is equivalent to $RM(d\cap X)=RM(d)$.
\bem The set of lines can be identified with $G/B\times G/B$, as $gBh^{-1}=g'Bh'^{-1}$ implies $gh^{-1}B^{h^{-1}}=g'h'^{-1}B^{h'^{-1}}$, and hence $h^{-1}h'\in N_G(B)=B$; similarly we obtain $g^{-1}g'\in N_G(B)=B$. It follows that the family of lines has Morley rank $2\,RM(G)-2\,RM(B)$.\ebem
\defn Let $X$ be a definable subset of $G$, and $g\in G$. Then $D_g(X)$ is the set of lines containing $g$ and generically contained in $X$.\edefn
\bem Clearly,
$$RM(X)\ge RM(\bigcup_{d\in D_g(X)}(X\cap d))=RM(D_g(X))+RM(B),$$
whence $RM(D_g(X))\le RM(X)-RM(B)$.\ebem

\defn Let $X$ be a definable subset of $G$. An element $g\in G$ is a {\em Fr\'econ element} for $X$ if $RM(D_g(X)) = RM(X) - RM(B)$. The set of Fr\'econ elements for $X$ is denoted by $\F(X)$.\edefn
In particular, if $RM(X)<RM(B)$ we have $\F(X)=\emptyset$.
\bem As a line is determined by two points, there are at most $2\,RM(X)-2\,RM(B)$ lines generically contained in $X$. Any such line contains at most $2\,RM(B)$ Fr\'econ points. Counting the pairs $(g,d)$ where $g\in\F(X)$ and $d\in D_g(X)$, we obtain
$$RM(\F(X))+RM(X)-RM(B)\le 2\,RM(X)-2\,RM(B)+RM(B),$$
whence $RM(\F(X))\le RM(X)$.\ebem

\bem If $X$ and $X'$ are definable subsets of $G$, then any line generically contained in $X\cup X'$ must be generically contained in $X$ or in $X'$ by connectedness of $B$. It follows that $\F(X\cup X')=\F(X)\cup\F(X')$.\ebem

\lmm\begin{enumerate}\item If $X\sim Y$, then $\F(X)=\F(Y)$.
\item If $g\in\F(X)$, then $hg\in\F(hX)$ and $gh\in\F(Xh)$.
\item We have $d\in D_g(X)$ if and only if $d^{-1}\in D_{g^{-1}}(X^{-1})$. In particular $\F(X)^{-1}=\F(X^{-1})$.
\item If $g\in\F(X)$ and $X$ has Morley degree $1$, then $X\sim\bigcup D_g(x)$. In particuliar $g^{-1}X\sim X^{-1}g$.\end{enumerate}\elmm
\bew\begin{enumerate}\item If $g\in\F(X)\setminus\F(Y)$, then 
$$RM(\bigcup_{d\in D_g(X)}((d\cap X)\setminus Y))=RM(D_g(X))+RM(B)=RM(X),$$
contradicting $X\sim Y$.
\item Obvious, as lines are preserved under translation.
\item Obvious, noting that for a line $d$ its inverse $d^{-1}$ is again a line.
\item As $g\in\F(X)$ we have 
$$RM(X)=RM(D_g(X))+RM(B)=RM(X\cap \bigcup D_g(X))\,;$$
since $X$ has Morley degree $1$ we obtain 
$$RM(X\setminus\bigcup D_g(X))<RM(X).$$
Conversely,
$$\begin{aligned}RM(\bigcup D_g(X)\setminus X)&=RM(\bigcup_{d\in D_g(X)}(d\setminus X))\\
&\le RM(D_g(X))+(RM(B)-1)<RM(X).\end{aligned}$$
Hence $X\sim\bigcup D_g(X)$.

If $g\in\F(X)$, then $1\in\F(g^{-1}X)$ and $1\in\F(X^{-1}g)$. So
$$g^{-1}X\sim\bigcup_{d\in D_1(g^{-1}X)}d=\bigcup_{d\in D_1(g^{-1}X)}d^{-1}=\bigcup_{d^{-1}\in D_1(X^{-1}g)}d^{-1}\sim X^{-1}g,$$
as the lines containing $1$ are subgroups.\qedhere\end{enumerate}\ebew
\satzli\label{p11} Let $G$ be a connected full Frobenius group of finite Morley rank, and $X$ a subset of $G$ which is not almost equal to $G$. Then either $\F(X)$ is contained in a finite union of translates of proper definable subgroups of $G$, or $G$ is not special and 
$$RM(\F(X))\le\frac{RM(G)+RM(B)}2-1<\frac34\,RM(G)-1.$$\esatzli
\bew Let $X$ be a definable subset of $G$ such that $\F(X)$ is not contained in a finite union of proper definable subgroups of $G$. We may suppose that $X$ is of Morley degree $1$; translating by $g^{-1}$ for some $g\in\F(X)$ we may furthermore assume that $1\in \F(X)$, and hence $X\sim X^{-1}$. Then for any $g\in\F(X)$ we have $g^{-1}X\sim X^{-1}g\sim Xg$, and the definable subgroup 
$$H=\{g\in G:\exists h\in G:gX\sim Xh\}$$
contains $\F(X)^{-1}=\F(X^{-1})=\F(X)$. So $H=G$ and $X$ is almost twistedly normal. By Proposition \ref{p7} there is a unique definable automorphism $\sigma$ such that $gX\sim X\sigma(g)$ for all $g\in G$, and $\sigma(g)=g^{-1}$ for all $g\in\F(X)$. So $\sigma^2$ fixes $\F(X)$ and hence $G$. It follows that $\sigma$ is involutive.

As $\sigma$ is not the identity, $G$ is not special by Proposition \ref{p6}, and $$RM(\F(X))\le RM(I)\le\frac{RM(G)+RM(B)}2-1<\frac34\,RM(G)-1.\qedhere$$\ebew

\bem In fact, Fr\'econ works with a set $X$ such that $X\sim\F(X)$. This gives some better bounds in Proposition \ref{p11}; we hope that separating $X$ and $\F(X)$ might allow a generalization of Fr\'econ's method. Note that the only currently known way to obtain a non-generic set $X$ with $\F(X)$ big, presented in the next section, automatically yields even $X\subseteq\F(X)$.\ebem

\section{Commutators}
In this section, we shall see that a simple full Frobenius group $G$ of finite Morley rank with abelian full complement $B$ contains a definable subset $X$ such that $RM(\F(X))=2\,RM(B)$, unless $RM(G)>2\,RM(B)+1$. This will quickly yield the main theorem.
\satzli\label{p13} Let $G$ be a simple full Frobenius group of finite Morley rank with an abelian full complement $B$. If $RM(G)=2\,RM(B)+1$, there is a definable subset $X$ with $RM(X)=2\,RM(B)$ and $X\subseteq\F(X)$.\esatzli
\bew Let $g=[h,h']$ be a non-trivial commutator, and put
$$X=\{x\in G:\exists\,y\in G\ [x,y]=g\}.$$
Note that 
$$\begin{array}{ll}
{[ux,y]=\ y^{-ux}y=\ \,y^{-x}y=[x,y]}&\mbox{if and only if }u\in C_G(y),\mbox{ and}\\
{[x,vy]=x^{-1}x^{vy}=x^{-1}x^y=[x,y]}&\mbox{if and only if }v\in C_G(x).\end{array}$$

Suppose first that $RM(X)=RM(G)$. For any $x\in X$ the set $\{y\in G:[x,y]=g\}=C_G(x)y$ has Morley rank $RM(B)$, and every non-trivial conjugacy class has Morley rank $RM(G)-RM(B)$. It follows that
$$\{(x,y)\in G\times G:[x,y]\mbox{ is conjugate to }g\}$$
has Morley rank 
$$RM(X)+RM(B)+RM(G)-RM(B)=2\,RM(G)\,;$$
it is thus generic in $G\times G$. It follows that for independent generic $x$ and $y$ both $[x,y]$ and $[y,x]$ are conjugate to $g$. So $[x,y]$ is conjugate to $[y,x]=[x,y]^{-1}$, a contradiction. Thus $RM(X)<RM(G)$.

Let $B$ be the Borel containing $h$, and for every $y\in Bh'$ let $B(y)$ be the Borel containing $y$. Then for every $x\in B(y)h$ we have
$$[x,y]=[h,y]=[h,h']=g,$$
whence $\bigcup_{y\in Bh'}B(y)h\subseteq X(g)$. Moreover, if $y'\in Bh'$ with $y\not=y'$, then $B(y)\not= B(y')$, as otherwise 
$$y'y^{-1}\in (Bh')(Bh')^{-1}=B,$$
so $B\cap B(y)$ contains the two points $1$ and $y'y^{-1}$. It follows that $h\in B=B(y)$, and $g=[h,y]=1$, a contradiction.

Thus $D_h(X)\supseteq\{B(y)h:y\in Bh'\}$ and 
$$RM(D_h(X))\ge RM(Bh')=RM(B).$$
Hence
$$RM(X)\ge RM(D_h(X))+RM(B)\ge 2\,RM(B),$$
and we must have equality. So $h\in\F(X)$; as $h\in X$ can be chosen arbitrarily, $X\subseteq\F(X)$.\ebew

\satz Let $G$ be a simple full Frobenius group of finite Morley rank and with abelian full complement $B$. Then $RM(G)> 2\,RM(B)+1$. In particular $RM(G)>3$.\esatz
\bew We already know that $RM(G)>2\,RM(B)$, so suppose $RM(G)=2\,RM(B)+1$. By Proposition \ref{p13} there is a definable set $X$ with $X\subseteq\F(X)$ and $RM(X)=2\,RM(B)$. Now
$$RM(\F(X))\ge 2\,RM(B)>\frac{RM(G)+RM(B)}2-1.$$
So $\F(X)$ and hence $X$ is contained in a finite union of translates of proper definable subgroups of $G$ by Proposition \ref{p11}. So there is a definable connected subgroup $H$ with $RM(H)=2\,RM(B)$ containing generically a line $d\in D_g(X)$ for some $g\in\F(X)$. But then $d^{-1}d$ is a Borel subgroup $B$ contained in $H$. Then $H$ is a full Frobenius group with full complement $B$ and of rank $2\,RM(B)$, a contradiction.\ebew

\section*{Acknowledgements}
The author should like to thank Bruno Poizat for many helpful discussions about Fr\'econ's proof.

\end{document}